\documentclass[11pt]{article}
\usepackage{amssymb,amsmath}
\newtheorem{thm}{Theorem}
\newtheorem{cor}[thm]{Corollary}
\newtheorem{prop}[thm]{Proposition}
\newtheorem{lemma}[thm]{Lemma}

\newcommand {\bC} {\mathbb {C}}
\newcommand {\bR} {\mathbb {R}}

\newcommand {\cD} {\mathcal {D}}

\newcommand {\cK} {\mathcal {K}}
\newcommand{\fT} {\mathfrak T}

\newcommand {\Tr} {\text {Tr}}

\def\mapright#1{\smash{\mathop{\longrightarrow}\limits^{#1}}}
\newcommand {\proof} {\noindent{\it Proof. }}

\begin{document}

\title{A note on noncommutative holomorphic and harmonic functions
on the unit disk.}

\author{\ \\S\l awomir Klimek \\ Department of Mathematics \\
Indiana University Purdue University Indianapolis\\ 402 N. Blackford St.
\\ Indianapolis, IN 46202 USA}
\date{\today }
\maketitle

\begin{abstract}
\vskip 0.3cm
\noindent We study noncommutative versions of holomorphic and harmonic 
functions on the unit disk.

\end{abstract}

\section{Introduction}\label{intesec}

The objective of this paper is to determine a complex structure on
the noncommutative disk $C(D_q)$, the q-deformation of the unit disk $D$. This noncommutative disk is a $C^*$-algebra that is a subalgebra of the quantum group $SU_q(2)$ and can be conveniently described using generators and (a quadratic) relation. It has been studied fairly extensively in the
literature - see \cite{KL}, \cite{NN}, \cite{SSV1}, \cite{SSV2}, \cite{SSV3}, \cite{SSV4}, \cite{SSV5}, and references therein. In particular the family of quantum disks $C(D_q)$ as $q$ varies, forms a deformation of the commutative disk, corresponding to $q=1$. To determine a complex structure on $C(D_q)$ we define and study partial derivatives on $C(D_q)$ and in particular the concept of holomorphic noncommutative functions.

The series of papers by Shklyarov and collaborators is very much related in spirit to our paper but is technically quite different, much more algebraic.
A similar study of a complex structure on the noncommutative plane is contained in \cite{RW}.

An important point of view of this paper is that we work in a concrete representation of $C(D_q)$ in a Hilbert space $H^2(D, d\mu)$ of holomorphic functions on the unit disk $D$, square integrable with respect to a certain measure. The algebra $C(D_q)$ is in this representation realized as the algebra of Toeplitz operators with continuous symbols. Also, we use this representation to realize different operations (scaling, derivatives, integral...) on $C(D_q)$ as coming from operators in $H^2(D, d\mu)$. This is well suited for operations that are only densely defined as it allows for good control over domains.

It turns out that there are two natural notions of holomorphic structure on the
quantum disk which we call weak and strong. Weakly holomorphic noncommutative functions directly correspond to ordinary holomorphic functions while the strongly holomorphic ones come from the scaled disk $\frac{1}{q}D$. We also study noncommutative harmonic functions. Just as ordinary two dimensional harmonic functions, their quantum counterparts on the unit disk can be written as a sum of holomorphic  and antiholomorphic part. They exhibit many of the familiar properties like a maximum principle.
The paper is organized as follows. In Section \ref{qudsec} we recall the definition of the quantum disk
and in particular we study in depth its representation using Toeplitz operators.
Section \ref{calcsec} contains the definition and our study of the properties of the derivatives 
and the integral on the quantum unit disk $C(D_q)$.
Finally in Section \ref{funsec} we introduce and study quantum holomorphic,
antiholomorphic and harmonic functions on the unit disk.

\section{Quantum unit disk}\label{qudsec}

In this section we review $C^*$-algebraic aspects of the quantum unit disk
$C(D_q)$. It is defined as the universal unital  $C^*$-algebra generated by
a generator $z$, and its conjugate denoted by $\bar z$, and satisfying the following
relation: $\bar zz=qz\bar z+(1-q)$. Symbolically:
\begin{eqnarray}
C(D_q):=<z,\bar z\ |\  \bar zz=qz\bar z+(1-q)>
\label{relationsref}
\end{eqnarray}

We will restrict ourself to $0\leq q<1$. Let us briefly recall the construction
of the universal $C^*$-algebra. If $a$ is a polynomial in $z,\bar z$ we define its norm
as the supremum of $||\rho(a)||$ over all Hilbert space representations $\rho$
satisfying the relation. One verifies that this defines a  sub-$C^*$-norm and the corresponding
completion mod the null space gives the universal $C^*$-algebra.

Notice that if $q=0$ the relation $\bar zz=1$ is the defining relation of the
standard  Toeplitz algebra $\fT$ - see \cite{F}. If $q=1$ the relation becomes the commutativity statement
$\bar zz=z\bar z$. Additionally, since $||\bar zz||=||z\bar z||=||z||^2$ we get
$||z||^2=q||z||^2+(1-q)$, which implies that $||z||=1$.
It is then natural to define $C(D_1)$ to be the algebra of continuous 
functions  $C(D)$ on the unit disk $D:=\{\zeta\in\bC:|\zeta|\leq 1\}$.

\begin{thm}\label{structhm}(see \cite{KL})

Let $\{e_n\}$ be the canonical basis in $l_2$, $n=0,1,2\ldots$, and let
$U:l_2\to l_2$ be the following weighted unilateral shift: $Ue_n=\sqrt{1-q^{n+1}}e_{n+1}$.
Then
$C(D_q)\cong C^*(U)$, where $C^*(U)$ is the $C^*$-algebra generated by $U$.

\end{thm}
\proof As noted above the case of $q=0$ is the standard Toeplitz algebra case.
If $q>0$
we calculate explicitly that $U^*e_{n+1}=\sqrt{1-q^{n+1}}\,e_{n}$, and consequently
$U^*Ue_n=(1-q^{n+1})\,e_{n}$ and $UU^*e_{n}=(1-q^{n})\,e_{n}$, which verifies
that $U$ gives a representation for $C(D_q)$. To verify that this indeed is the defining
representation one needs to classify all irreducible representations.
This was done in \cite{KL}.
$\square$

\begin{cor}(see \cite{KL})

We have the exact sequence:
\begin{eqnarray}
0\ \mapright{}\ \cK
\ \mapright{}\ C(D_q)\ \mapright{\sigma}\ C(\partial D)\ \mapright{}\ 0,
\label{exseq}
\end{eqnarray}
where $\cK$ is the ideal of compact operators in $l_2$ and $\partial D$, the boundary
of $D$, is the unit circle.
\end{cor}
\proof The corollary follows from the general theory of weighted shifts \cite{C}.
Briefly, the commutator $[U,U^*]$ is compact since its eigenvalues
$(1-q)q^n\to 0$ as $n\to\infty$, and, since the $C^*$-algebra is irreducible, it contains all compact
operators. The quotient $C(D_q)/\cK$ is generated by the unitary operator $[U]$
the spectrum of which is the full unit circle. For details see \cite{KL}.
$\square$

In the exact sequence (\ref{exseq}) the map 
$\sigma:C(D_q)\to C(\partial D)$ is called the symbol map.

\begin{prop}
{\ }

The $C^*$-algebras $C(D_q)$ are isomorphic to each other and
for every $q$, $0\leq q<1$, we have
$C(D_q)\cong \fT$, where $\fT$ is the  Toeplitz algebra. 
\end{prop}
\proof Let $V:l_2\to l_2$ be the unilateral shift $Ve_n=e_{n+1}$, so that
$C^*(V)=\fT$. Notice that $U-V$ is a weighted shift $(U-V)e_n=\lambda_ne_{n+1}$,
where weights $\lambda_n\to 0$ as $n\to\infty$. Consequently $U-V$ is
a compact operator. The proposition now follows from the exact sequence (\ref{exseq}).
$\square$

The $C^*$-algebras $C(D_q)$ are also continuous in $q$ in the following sense:

\begin{thm}(see \cite{NN})

The $C^*$-algebras $C(D_q)$ for $0\leq q<1$ and $C(D)$ for $q=1$ form a continuous field of $C^*$-algebras with the space of cross-sections obtained by completing the space of polynomials in $z$ and $\bar z$ with coefficients
which are continuous functions of $q$. 
\end {thm}
\proof This was done by Nagy and Nica in \cite{NN} for an even bigger
range $-1\leq q\leq 1$.
$\square$

As the final part of this section we will discuss another useful representation
of $C(D_q)$. It is using Toeplitz operators and is implicitly contained
in \cite{KL} but we work out the details here. 

Consider the following measure on the unit disk $D$:
\begin{eqnarray}
d\mu(\zeta)=\prod\limits_{i\geq 0}(1-|\zeta|^2q^{i+1})\sum\limits_{m\geq 0}q^m
\delta_{|\zeta|^2=q^m}(\zeta),
\end{eqnarray}
where $\delta_{|\zeta|^2=r^2}$ is the normalized Lebesgue measure
on the circle $|\zeta|^2=r^2$. Let $H^2(D, d\mu)\subset L^2(D, d\mu)$
be the closed subspace consisting of holomorphic functions and
let $P$ be the corresponding orthogonal projection. If $f\in C(D)$ we define the Toeplitz
operator $T(f):H^2(D, d\mu)\to H^2(D, d\mu)$ with symbol $f$, by:
$T(f)=PM(f)P$, where $M(f)$ is the multiplication by $f$.

\begin{thm}\label{holorep}(see \cite{KL})  

With the above notation, the $C^*$-algebra generated by 
$\{T(f)\}$, $f\in C(D)$, is naturally isomorphic with
$C(D_q)$. The isomorphism is
determined by identifications:
$z=T(\zeta)$, $\bar z=T(\bar\zeta)$
\end {thm}
\proof It follows from the definition that
$||T(f)||\leq{\mathop{\sup}_{\zeta\in D}}|f(\zeta)|$. Since polynomials in
$\zeta$ and $\bar\zeta$ are dense in $C(D)$ and $T(\bar\zeta^m\zeta^n)=
(T(\zeta)^*)^m\,(T(\zeta))^n$, we see that the algebra generated by Toeplitz operators
is in fact generated by the single operator $T(\zeta)$. Next, because $d\mu$
is rotationally invariant, the functions $\zeta^n$ are mutually orthogonal
and form an unnormalized basis in $H^2(D, d\mu)$. To find an orthonormal basis we compute
\begin{eqnarray*}
\int|\zeta|^{2n}\,d\mu(\zeta)&=&\sum_{m\geq 0}q^mq^{nm}\prod_{i\geq 0}(1-q^mq^{i+1})=\\
&=&\prod_{i\geq 0}(1-q^{i+1})\left(1+\sum_{m\geq 1}\frac{q^{m(n+1)}}
{\prod_{1\leq k\leq m}(1-q^k)}\right).
\end{eqnarray*}
Using the Euler's identity:
\begin{eqnarray}
\frac{1}{\prod_{i\geq 0}(1-xq^i)}=1+\sum_{m\geq 1}\frac{x^m}
{\prod_{1\leq k\leq m}(1-q^k)},\label{euler}
\end{eqnarray}
gives $\int|\zeta|^{2n}\,d\mu(\zeta)=1$ if $n=0$, and, for $n\geq 1$:
\begin{eqnarray*}
\int|\zeta|^{2n}\,d\mu(\zeta)=\frac{\prod_{i\geq 0}(1-q^{i+1})}
{\prod_{i\geq 0}(1-q^{n+i+1})}=\prod_{i= 0}^{n-1}(1-q^{i+1}).
\end{eqnarray*}
It follows that the measure $d\mu$ is probabilistic and the following is
an orthonormal basis in $H^2(D, d\mu)$:
\begin{eqnarray}
e_n=\begin{cases} 
1& \text{if $n=0$}\\
\frac{\zeta^n}{\sqrt{\prod_{i=0}^{n-1}(1-q^{i+1})}}& \text{if $n\geq 1$}
\end{cases}\label{basis}
\end{eqnarray}
We now find the matrix elements of $T(\zeta)$ with respect to the basis $e_n$:
\begin{eqnarray*}
T(\zeta)e_n&=&\zeta e_n=\frac{\zeta^n}{\sqrt{\prod_{i=0}^{n-1}(1-q^{i+1})}}=
\sqrt{1-q^{n+1}}\frac{\zeta^n}{\sqrt{\prod_{i=0}^{n}(1-q^{i+1})}}=\\
&=&\sqrt{1-q^{n+1}}\,e_{n+1}.
\end{eqnarray*}
So the matrix elements of $T(\zeta)$ are equal to that of $U$ of the structure theorem \ref{structhm},
which concludes the proof.
$\square$

From now on we will identify $C(D_q)$ with the concrete algebra generated by
Toeplitz operators in $H^2(D, d\mu)\subset L^2(D, d\mu)$.  For future reference
we recall here the definition of the Bergman kernel $K(\zeta,\bar\eta)$ for $H^2(D, d\mu)$. 
It is the integral kernel of the projection $P$ so it has the reproducing property:
\begin{eqnarray}
\int K(\zeta,\bar\eta)\phi(\eta)\,d\mu(\eta)=\phi(\zeta),\label{reprodprop}
\end{eqnarray}
where $\phi(\zeta)\in H^2(D, d\mu)$.
It can be explicitly computed using a basis in $H^2(D, d\mu)$, 
for example the one given by (\ref{basis}). We obtain:
\begin{eqnarray}
K(\zeta,\bar\eta)&=&\sum_{n=0}^\infty e_n(\zeta)\overline{e_n(\eta)}=
1+\sum_{n\geq 1}\frac{(\zeta\bar\eta)^n}{\prod_{1\leq k\leq n}(1-q^k)}=\nonumber\\
&=&=\frac{1}{\prod_{i\geq 0}(1-\zeta\bar\eta q^i)}.
\label{kernel}
\end{eqnarray}
In the above we again used the Euler identity (\ref{euler}).

By construction, the space
of polynomials in Toeplitz operators is dense in $C(D_q)$. More is actually true
as spelled out in the next statement.

\begin{prop}\label{densityprop}
{\ }

The subspace of Toeplitz operators $T(f),\ f\in C(D)$ is dense in $C(D_q)$.
\end {prop}
\proof It follows from the defining relation of $C(D_q)$ that the linear span
of $\bar z^mz^n$, $m,n\geq 0$, forms a dense subalgebra of $C(D_q)$. Indeed,
since $z\bar z$ expresses linearly in terms of $\bar zz$ we can rearrange any 
polynomial in $z,\bar z$ so that powers of $\bar z$ come first.
But $\bar z^mz^n=T(\bar\zeta^m\zeta^n)$ and the claim follows.
$\square$

\section{Calculus on $C(D_q)$}\label{calcsec}

In this section we introduce calculus on the quantum unit disk. In the following we assume that $q>0$. Formal aspects of the calculus on $C(D_q)$ can be found in
\cite{CHZ} as well as in \cite{SSV1,SSV2,SSV3,SSV4,SSV5}. We concentrate here on
issues of domains for various unbounded operators and we will always identify $C(D_q)$ with
the concrete algebra of Toeplitz operators of Theorem \ref{holorep}.

Let $\cD:=\{\phi\in H^2(D, d\mu):\ \phi(\zeta/q)\in H^2(D, d\mu)\}$.
Clearly $\cD$ is a dense subspace in $H^2(D, d\mu)$ containing all polynomials,
or more generally entire functions. We define a scaling operator 
$j: H^2(D, d\mu)\to H^2(D, d\mu)$ by the formula:
\begin{eqnarray*}
j\phi(\zeta):=\phi(q\zeta).\label{jformula}
\end{eqnarray*}
The operator $j$ is bounded, one-to-one and $\text {Ran}j=\cD$.
Using the defining formula (\ref{basis}) we have 
\begin{eqnarray}
je_n=q^ne_n,\label{jbasis}
\end{eqnarray}
so that $j$ is a self-adjoint
compact operator. Since $z\bar ze_n=(1-q^n)e_n$ and $\bar zze_n=(1-q^{n+1})e_n$ we have
\begin{eqnarray}
z\bar z=1-j,\ \bar zz=1-qj.\label{jzrelation}
\end{eqnarray}

An element $a\in C(D_q)$ is called {\it scalable} if $J(a):=j^{-1}aj$ is a bounded operator.
We have a simple proposition:
\begin{prop}\label{dinvprop} 
{\ }

The operator $j^{-1}aj$ is bounded iff $a$ maps $\cD$ to $\cD$.
\end {prop}
\proof  If $a$ preserves $\cD$ then $j^{-1}aj$ is defined everywhere. To show that it is bounded,
we use the closed graph theorem which implies that we need to verify that if $x_n\to x$
and $y_n:=J(a)x_n\to y$ then $J(a)x=y$. Since $j$ is continuous we have $jy_n=ajx_n\to jy$.
But $aj$ is continuous so $ajx_n\to x$ and consequently $ajx=jy$, which is what we wanted.
The converse statement is straightforward.
$\square$

\medskip
We write $C_s(D_q)$ for the set of scalable elements of $C(D_q)$. The proposition below
shows that $C_s(D_q)$ is a subalgebra of $C(D_q)$ containing $\text {Pol}(D_q)$, the algebra
of polynomials in $z,\bar z$. However, examples below show that $C_s(D_q)$ is not closed with respect to taking adjoints and inverses.

\begin{prop}\label{Jprop}
{\ }

With the above notation we have: 
\begin{eqnarray}
J(1)=1,\ J(\bar z)=q\bar z,\ J(z)=q^{-1}z.\label{jproperties}
\end{eqnarray}
If $a,b\in C_s(D_q)$ then $ab\in C_s(D_q)$ and $J(ab)=J(a)J(b)$.
\end {prop}
\proof  The proof consist of straightforward computations verifying each of the properties.
For this we need explicit formulas for $z, \bar z$. Theorem \ref{holorep} implies
\begin{eqnarray}
z\phi(\zeta)=\zeta\phi(\zeta),\label{zformula}
\end{eqnarray}
while the structure Theorem \ref{structhm} says that
\begin{eqnarray}
ze_n=\sqrt{1-q^{n+1}}\,e_{n+1},\label{zbasis}
\end{eqnarray}
where $e_n$ were defined in (\ref{basis}). Taking the adjoint gives
\begin{eqnarray}
\bar ze_n=\sqrt{1-q^{n}}\,e_{n-1},\label{zbarbasis}
\end{eqnarray}
(the right-hand side is defined to be 0 when $n=0$). This implies that 
$\bar z\zeta^n =(1-q^n)\zeta^{n-1}$, which in turn gives:
\begin{eqnarray}
\bar z\phi(\zeta)=\frac{\phi(\zeta)-\phi(q\zeta)}{\zeta}.\label{zbarformula}
\end{eqnarray}
A sample calculation verifying one of the statements of the proposition follows:
\begin{eqnarray*}
J(z)\phi(\zeta)&=&j^{-1}zj\phi(\zeta)=zj\phi(\zeta/q)=q^{-1}\zeta j\phi(\zeta/q)\\
&=&q^{-1}\zeta\phi(\zeta)=q^{-1}z\phi(\zeta)
\end{eqnarray*}

$\square$

We are now going to look at examples to illustrate some subtleties of the notion of scalability.
First notice that $1-q\bar z$ is invertible since $||q\bar z||=q<1$. The inverse $a:=(1-q\bar z)^{-1}$
is clearly in $C(D_q)$ and is scalable because $J(a)=(1-q^2\bar z)^{-1}$. However, 
$a^* =(1-qz)^{-1}$ is not scalable as $J(a^*)=(1-z)^{-1}$ is unbounded. Next consider
$b:=1-qz$. Clearly $b\in C(D_q)$, $b$ is scalable, $b$ is invertible, and the inverse of $b$
is in $C(D_q)$. But since $b^{-1}=a^*$, $b^{-1}$ is not scalable.

Next we introduce two operators $\delta, \bar\delta$ in $H^2(D,d\mu)$ that will be used to define Dolbeault - type operators $\partial,\bar\partial$ on $C(D_q)$. The precise form of $\delta, \bar\delta$ is dictated by the desired properties of $\partial,\bar\partial$ as described in the Proposition
\ref{deltaprop}.

The operators $\delta, \bar\delta$ are defined to be unbounded operators in $H^2(D,d\mu)$ with domains both equal to $\cD$ and given by the following formulas using $z,\bar z,j$:
\begin{eqnarray}
\bar\delta={(q-1)}^{-1}j^{-1}z={(q-1)}^{-1}q^{-1}zj^{-1},\label{deltadef}
\end{eqnarray}
\begin{eqnarray}
\delta={(1-q)}^{-1}j^{-1}\bar z={(1-q)}^{-1}q\bar zj^{-1}.\label{bardeltadef}
\end{eqnarray}

\begin{prop}
{\ }

With the above notation, the operators $\delta, \bar\delta$ are closed
(on $\cD$).
\end{prop}
\proof To show that $\bar\delta$ is closed the following needs to be demonstrated:
if $\phi_n\to \phi$, $\phi_n\in\cD$ and $\psi_n:=j^{-1}z\phi_n\to \psi$,
then $\phi\in\cD$ and $j^{-1}z\phi=\psi$. Applying $j$ to $\psi_n$ and
using the continuity of $j$ gives $z\phi_n\to j\psi$. On the other hand, since
$z$ is continuous, we have $z\phi_n\to z\phi$. Consequently $z\phi=j\psi$,
which means that $z\phi\in\cD$ and $j^{-1}z\phi=\psi$. What's left is
to show that $\phi\in\cD$. Applying $\bar z$ to both sides of $z\phi=j\psi$
and using (\ref{jzrelation}) and (\ref{jproperties}) we obtain
\begin{eqnarray*}
\phi=qj(1-qj)^{-1}\bar z\psi,
\end{eqnarray*}
which concludes the proof that $\bar\delta$ is closed. The proof for
$\delta$ is analogous with the exception of the fact that $z\bar z$ has a kernel. The analog of the above formula works on the orthogonal complement
of that kernel. The proof is then concluded by observing that the
kernel of $z\bar z$ is one dimensional and is contained in $\cD$.
$\square$

Using the equations (\ref{jformula}), (\ref{zformula}) and (\ref{zbarformula}), we obtain the following explicit descriptions of the operators 
$\delta, \bar\delta$:
\begin{eqnarray}
\bar\delta\phi(\zeta)={(q-1)}^{-1}\zeta/q\,\phi(\zeta/q)
\end{eqnarray}
\begin{eqnarray}
\delta\phi(\zeta)=\frac{\phi(\zeta)-\phi(\zeta/q)}{(1-1/q)\zeta}.
\end{eqnarray}
Optionally, when working with the operators $\delta, \bar\delta$ one can use their matrix elements,
obtained using  (\ref{jbasis}), (\ref{zbasis}) and (\ref{zbarbasis}):
\begin{eqnarray*}
\bar\delta e_n={(q-1)}^{-1}q^{-(n+1)}\sqrt{1-q^{n+1}}\,e_{n+1},\
\end{eqnarray*}
\begin{eqnarray*}
\delta e_n=\frac{1-1/q^n}{1-1/q}\sqrt{1-q^{n}}\,e_{n-1}.
\end{eqnarray*}

Now we use those operators to define complex structure on $C(D_q)$ - for this we need the analogs of the usual complex derivatives $\partial,\bar\partial$. They are defined using scaled commutators with $\delta, \bar\delta$ as follows. If $a$ is scalable we define $\bar\partial(a)$, $\partial(a)$ to be (in general unbounded) linear operators defined on $\cD$ by:
\begin{eqnarray*}
\bar\partial(a)=\bar\delta a-J(a)\bar\delta,
\end{eqnarray*}
\begin{eqnarray*}
\partial(a)=\delta a-J(a)\delta.
\end{eqnarray*}
Proposition \ref{dinvprop} assures that $\bar\partial(a)$, $\partial(a)$ are well defined operators on $\cD$. For general $a$ i.e. not necessarily
scalable, $\bar\partial(a)$, $\partial(a)$ make sense only as quadratic forms
- see below. We will use those quadratic forms in the discussion of quantum holomorphic and harmonic functions in the next section.

The following proposition summarizes the main properties of the operators $\partial,\bar\partial$.

\begin{prop}\label{deltaprop} 
{\ }

With the above notation we have, assuming $a,b\in C_s(D_q)$: 
\begin{eqnarray}
\bar\partial(1)=0,\ \bar\partial(\bar z)=1,\ \bar\partial(z)=0,\ 
\bar\partial(ab)=(\bar\partial a)b+J(a)(\bar\partial b),
\end{eqnarray}
\begin{eqnarray}
\partial(1)=0,\ \partial(\bar z)=0,\ \partial(z)=1,\ 
\partial(ab)=(\partial a)b+J(a)(\partial b).
\end{eqnarray}
In particular if $a\in \text {Pol}(D_q)$ then 
$\partial a,\bar\partial a\in \text {Pol}(D_q)$.
\end {prop}
\proof  The proof again consist of straightforward verifications using definitions. Below we show calculations of the action of the operators $\partial,\bar\partial$ on $z,\bar z$ that utilize commutation relations among $j,z,\bar z$. All the manipulations with unbounded operators make
sense pointwise on $\cD$.
\begin{eqnarray*}
\partial(z)&=&\delta z-J(z)\delta={(1-q)}^{-1}(j^{-1}\bar z z-q^{-1}zj^{-1}\bar z)=\\
&=&{(1-q)}^{-1}j^{-1}(\bar z z-z\bar z)={(1-q)}^{-1}j^{-1}(1-qj-1+j)=1\\
\end{eqnarray*}
\begin{eqnarray*}
\partial(\bar z)&=&\delta\bar z-J(\bar z)\delta={(1-q)}^{-1}(j^{-1}\bar z^2-q\bar zj^{-1}\bar z)=\\
&=&{(1-q)}^{-1}j^{-1}(\bar z^2-\bar z^2)=0\\
\end{eqnarray*}
The other two calculations are very similar.
$\square$

For future reference we note the formulas for action of $\partial,\bar\partial$
on monomials:
\begin{eqnarray}
\partial(\bar z^{n}z^{m})=q^{n-m+1}[m]_q\bar z^{n}z^{m-1}\label{partialmono}
\end{eqnarray}
and similar:
\begin{eqnarray}
\bar\partial(\bar z^{n}z^{m})=[n]_q\bar z^{n-1}z^{m}.\label{barpartialmono}
\end{eqnarray}
Here, and later in the paper, we use the notation
\begin{eqnarray*}
[n]_q:=\frac{1-q^n}{1-q}.
\end{eqnarray*}
Another set of useful formulas follows directly from the definitions:
\begin{eqnarray}
\partial a={(1-q)}^{-1}j^{-1}[\bar z,a]\label{paformula}
\end{eqnarray}
\begin{eqnarray}
\bar\partial a={(q-1)}^{-1}j^{-1}[z,a]\label{bpaformula}
\end{eqnarray}
For example, it follows from those formulas that $\partial a,\bar\partial a$ are closable since the domains of the adjoints clearly contain $\cD$ and so are dense. Another application of (\ref{paformula}) and (\ref{bpaformula})
is in the following definition of the derivatives $\partial a,\bar\partial a$ as quadratic forms for a general, not necessarily scalable $a\in C(D_q)$. They are defined on $\cD$ as
\begin{eqnarray}
Q_{\partial a}(\phi):={(1-q)}^{-1}(j^{-1}\phi,[\bar z,a]\phi),\label{Qpa}
\end{eqnarray}
\begin{eqnarray}
Q_{\bar\partial a}(\phi):={(q-1)}^{-1}(j^{-1}\phi,[z,a]\phi).\label{Qbpa}
\end{eqnarray}

We now turn to the definition and properties of the laplacian on $C(D_q)$.
There are two natural choices that we will look at using the formulas above:
\begin{eqnarray*}
\bar\partial\partial(\bar z^{n}z^{m})=q^{n-m+1}[m]_q[n]_q\bar z^{n-1}z^{m-1}.
\end{eqnarray*}
Similarly we obtain
\begin{eqnarray*}
\partial\bar\partial(\bar z^{n}z^{m})=q^{n-m}[m]_q[n]_q\bar z^{n-1}z^{m-1}.
\end{eqnarray*}
It follows that, at least on $\text{Pol}(D_q)$,
\begin{eqnarray}
\bar\partial\partial=q\partial\bar\partial.\label{twolaplace}
\end{eqnarray}
To define $\bar\partial\partial$ and $\partial\bar\partial$ for a larger class of elements
of $C(D_q)$ we proceed similarly to the way we defined $\partial,\bar\partial$. Let
$\cD_2:=\{\phi\in H^2(D, d\mu):\ \phi(\zeta/q^2)\in H^2(D, d\mu)\}$. Clearly  
$\cD_2=\text {Ran}j^2$, $\cD_2$ is dense and $\cD_2\subset\cD$. Also, just as in Proposition \ref{dinvprop}, if $a$ is scalable and $J(a)$ is scalable then $a$ maps $\cD_2$ into $\cD_2$. In particular,
$z,\bar z$ preserve $\cD_2$. Consequently, it follows from (\ref{deltadef}) and (\ref{bardeltadef})
that $\delta,\bar\delta:\cD_2\to\cD$. Thus if both $a$ and $J(a)$ are scalable then
$\bar\partial\partial(a)$ and $\partial\bar\partial(a)$ make sense as operators on $\cD_2$.
It can be easily verified that (\ref{twolaplace}) holds in this more general context i.e.
if $a,J(a)$ are scalable and $\phi\in\cD_2$ then
\begin{eqnarray*}
\bar\partial\partial(a)\phi=q\partial\bar\partial(a)\phi.
\end{eqnarray*}
In particular, the two laplacians have the same kernels and we will use whatever is more convenient when defining harmonic functions as extended kernels in the next section.

The last item in this section is integration on the quantum unit disk.
We define the integral $\int_{D_q}: C(D_q)\to\bR$ by
\begin{eqnarray}
\int_{D_q}a=\frac{\Tr(aj)}{\Tr(j)}
\end{eqnarray}
Using (\ref{jbasis}) we compute:
\begin{eqnarray*}
\Tr(j)=\sum_{n=0}^\infty q^n=\frac{1}{1-q}
\end{eqnarray*}
\begin{prop}
{\ }

$\int_{D_q}$ is a faithful state on $C(D_q)$ and
\begin{itemize}
\item $\int_{D_q}ab=\int_{D_q}J(b)a$
\item $\int_{D_q}J(a)=\int_{D_q}a$
\end{itemize}
\end{prop}
\proof Everything follows easily from the definitions.
$\square$

The integral is easy to work with as demonstrated in the following computation of its value 
on monomials. 
\begin{lemma}\label{intlemma}
{\ }

With the above notation we have
\begin{eqnarray*}
\int_{D_q}\bar z^nz^m=\delta_{n,m}\ \frac{1}{[n+1]_q}
\end{eqnarray*}
\end{lemma}
\proof
Using the canonical basis in $H^2(D, d\mu)$ we have
\begin{eqnarray*}
\int_{D_q}a=(1-q)\sum_{k=0}^\infty q^k(e_k,ae_k).
\end{eqnarray*}
It follows that $\int_{D_q}\bar z^nz^m=0$ if $n\ne m$. Using (\ref{zbasis}) we compute
\begin{eqnarray*}
\int_{D_q}\bar z^nz^n=(1-q)\sum_{k=0}^\infty q^k(e_k,\bar z^nz^ne_k)=
(1-q)\sum_{k=0}^\infty q^k||z^ne_k||=\\
=(1-q)\sum_{k=0}^\infty q^k(1-q^{k+1})(1-q^{k+2})\ldots(1-q^{k+n})=
\int_0^1f(y)\,d_qy,\\
\end{eqnarray*}
where $f(y)=(1-qy)(1-q^2y)\ldots(1-q^ny)$.
Here we used the Jackson's integral for a continuous function $f$ which is defined by:
\begin{eqnarray*}
\int_0^1f(y)\,d_qy:=(1-q)\sum_{k=0}^\infty q^k f(q^k)
\end{eqnarray*}
It has the property
\begin{eqnarray}
\int_0^1\delta_qg(y)\,d_qy=g(1)-g(0),\label{qintder}
\end{eqnarray}
where 
\begin{eqnarray}
\delta_qg(y)=\frac{g(y)-g(qy)}{y-qy}.\label{qder}
\end{eqnarray}
We use this property in our calculation. For
$g(y)=(1-y)(1-qy)\ldots(1-q^ny)$ we compute:
\begin{eqnarray*}
&\delta_qg(y)
=\frac{(1-y)(1-qy)\ldots(1-q^ny)-(1-qy)(1-q^2y)\ldots(1-q^{n+1}y)}{y(1-q)}=\\
&=(1-qy)(1-q^2y)\ldots(1-q^ny)\frac{1-y-1+q^{n+1}y}{y(1-q)}=-[n+1]_qf(y).\\
\end{eqnarray*}
It follows that
\begin{eqnarray*}
\int_0^1f(y)\,d_qy=\int_0^1\delta_q\left(\frac{-g(y)}{[n+1]_q}\right)\,d_qy=\frac{g(0)-g(1)}{[n+1]_q},
\end{eqnarray*}
which finishes the proof.
$\square$

The integral $\int_{D_q}$ and the derivatives $\bar\partial,\partial$ are tightly connected,
just as in the commutative case. This is illustrated by the following theorem - compare
also \cite{SSV1}.
\begin{thm} (Green's theorem)

If $a\in\text {Pol}(D_q)$ then
\begin{eqnarray*}
\int_{D_q}(\bar\partial a)=\frac{1}{2\pi i}
\int_{\partial D}\sigma(a)(\zeta)d\zeta.
\end{eqnarray*}
Here $\sigma:C(D_q)\to C(S^1)$ is the symbol map.
\end{thm}

\proof
It is enough to consider monomials of the following form:
\begin{eqnarray*}
a=\bar z^{n+1}z^n.
\end{eqnarray*}
Then, $\sigma(a)=\bar\zeta^{n+1}\zeta^n=\bar\zeta=e^{-i\theta}$, and
\begin{eqnarray*}
\frac{1}{2\pi i}\int_{\partial D}\sigma(a)(\zeta)\,d\zeta=
\frac{1}{2\pi i}\int_{0}^{2\pi}e^{-i\theta}\,d(e^{i\theta})=1.
\end{eqnarray*}
On the other hand $\partial(\bar\zeta^{n+1}\zeta^n)=[n+1]_q\bar\zeta^{n}\zeta^n$ 
by (\ref{barpartialmono}) and
\begin{eqnarray*}
\int_{D_q}(\bar\partial a)=[n+1]_q\int_{D_q}\bar\zeta^{n}\zeta^n=1
\end{eqnarray*}
by Lemma \ref{intlemma}.
$\square$

\section{Quantum holomorphic and harmonic functions}\label{funsec}

In this section we define quantum holomorphic and harmonic functions on the quantum unit disk $C(D_q)$.
We start with the following definition. An element
$a\in C(D_q)$ is called strongly holomorphic if $a$ is scalable and $\bar\partial a=0$. Similarly $a\in C(D_q)$ is called weakly holomorphic if $Q_{\bar\partial a}(\phi)=0$ for all $\phi\in\cD$, where the quadratic form $Q_{\bar\partial a}$ was defined in (\ref{Qbpa}). In the later definition we do not need to assume scalability of $a$. We denote  by $\text {Hol}(D_q)$ the space of weakly holomorphic elements of $C(D_q)$. We have the following simple proposition:
\begin{prop}\label{holofprop}
{\ }
\begin{itemize}
\item $a\in C(D_q)$ is weakly holomorphic iff $[z,a]=0$.
\item $a\in C(D_q)$ is strongly holomorphic iff $a$ is scalable and weakly holomorphic.
\end{itemize}
\end{prop}
\proof The formula (\ref{Qbpa}) and polarization imply the first part of the proposition. The second part is just a rephrasing of the definition.
$\square$

There are analogous definitions of antiholomorphic functions. An element
$a\in C(D_q)$ is called strongly antiholomorphic if $a$ is scalable and $\partial a=0$.
Similarly $a\in C(D_q)$ is called weakly antiholomorphic if $Q_{\partial a}(\phi)=0$ 
for all $\phi\in\cD$, where the quadratic form $Q_{\partial a}$ was defined in (\ref{Qpa}). Because $z$ and $\bar z$ scale differently, the following
analog of Proposition \ref{holofprop} looks a little different.

\begin{prop}\label{aholofprop}
{\ }

\begin{itemize}
\item $a\in C(D_q)$ is weakly antiholomorphic iff $[\bar z,a]=0$.
\item $a\in C(D_q)$ is strongly antiholomorphic iff $a$ is  weakly antiholomorphic.
\end{itemize}
\end{prop}
\proof The formula (\ref{Qpa}) and polarization imply the first part of the proposition. The second part is proved in the theorem below.
$\square$

The following is the main result describing holomorphic and antiholomorphic 
functions on the quantum unit disk. Notice that there is slight asymmetry between the notions of strongly holomorphic and antiholomorphic functions which disappears when $q=1$.
\begin{thm}\label{qholothm}
{\ }
\begin{enumerate}
\item  If $f\in C(D)$ is holomorphic inside $D$ then the corresponding Toeplitz operator $T(f)\in C(D_q)$ is weakly holomorphic.
\item  If $a\in C(D_q)$ is weakly holomorphic then there exist $f\in C(D)$ which is holomorphic inside $D$ such that $a=T(f)$.
\item  (maximum principle) If $a\in\text {Hol}(D_q)$  then 
$||a||_{D_q}=||\sigma(a)||_{\partial D}$
\item  The space $\text {Hol}(D_q)\subset C(D_q)$ is a Banach subalgebra isomorphic to the algebra $\text {Hol}(D)\subset C(D)$ of continuous functions on $D$ and holomorphic inside $D$.
\item The above statements are also true when the word holomorphic is replaced by antiholomorphic throughout.
\item If $a\in C(D_q)$ then $a$ is strongly antiholomorphic iff $a$ is  weakly antiholomorphic.
\end{enumerate}
\end{thm}
\proof We proof all items in order stated in the theorem:

1. This follows from Proposition \ref{holofprop} since if $f\in C(D)$ is holomorphic then:
\begin{eqnarray*}
zT(f)\phi(\zeta)=\zeta f(\zeta)\phi(\zeta)=T(f)z\phi(\zeta).
\end{eqnarray*}

2. For a weakly holomorphic $a\in C(D_q)$ we set $f(\zeta):=a\cdot 1(\zeta)\in H^2(D, d\mu)$.
In particular $f$ is holomorphic inside the disk $D$. 
Because $[z,a]=0$, we have inductively $a\zeta^n=f(\zeta)\zeta^n$,
so $a$ is equal to the Toeplitz operator $T(f)$ on the dense domain and consequently everywhere.
To obtain more information about $f$ we prove the following estimate:
\begin{eqnarray}
\sup_{\zeta\in D}|f(\zeta)|\leq ||T(f)||.\label{normestim}
\end{eqnarray}
To do it we consider the family of functions:
\begin{eqnarray*}
\phi_\eta(\zeta):=\frac{K(\bar\eta,\zeta)}{(K(\bar\eta,\eta))^{1/2}},
\end{eqnarray*}
where $K(\zeta,\bar\eta)$ is the reproducing kernel (\ref{kernel}).
It is easily seen that the functions $\phi_\eta(\zeta)$ belong to $H^2(D, d\mu)$ and have norm 1.
Using the reproducing property (\ref{reprodprop}) we compute:
\begin{eqnarray*}
(\phi_\eta,T(f)\phi_\eta)&=&
(K(\bar\eta,\eta))^{-1}\int K(\eta,\bar\zeta)K(\zeta,\bar\eta)f(\zeta)\,d\mu(\zeta)=\\
&=&(K(\bar\eta,\eta))^{-1}K(\bar\eta,\eta)f(\eta)=f(\eta).\\
\end{eqnarray*}
It follows that
\begin{eqnarray*}
\sup_{\zeta\in D}|f(\zeta)|=\sup_{\zeta\in D}|(\phi_\zeta,T(f)\phi_\zeta)|\leq ||T(f)||<\infty,
\end{eqnarray*}
so $f$ is bounded on $D$ and holomorphic inside it. But $a=T(f)$ belongs
to $C(D_q)$ so it is a limit of polynomials which implies that
$f\in C(D)$ as claimed.

3. It follows from the definition that $||T(f)||\leq{\mathop{\sup}_{\zeta\in D}}|f(\zeta)|$. 
On the other hand if $a=T(f)\in\text {Hol}(D_q)$ the estimate (\ref{normestim}) is valid and
so $||a||=||T(f)||=\sup_{\zeta\in D}|f(\zeta)|$. But $\sigma(T(f))=f|_{\partial D}$ and the
supremum of $|f(\zeta)|$ is achieved on the boundary ${\partial D}$ and the statement follows.

4. This is just a rephrasing of the previous items.

5. Notice that Propositions \ref{holofprop} and \ref{aholofprop} imply that if $a$ is weakly antiholomorpic then $a^*$ is weakly holomorphic, so all the statements follow by conjugation.

6. The previous item implies that if $a$ is weakly antiholomorphic then $a=T(f)$
where $f$ is antiholomorphic. For such $f$ we have $J(T(f(\zeta)))=T(f(q\zeta))$, so $T(f)$ is automatically scalable. Consequently, $a$ is strongly antiholomorphic.
$\square$

Notice that if $f$ is holomorphic then $J(T(f(\zeta)))=T(f(\frac{1}{q}\zeta))$ 
so that if $T(f)$ is strongly holomorphic then $f$ extends to a holomorphic function inside the disk 
$\frac{1}{q}D=\{\zeta\in\bC:|\zeta|\leq\frac{1}{q}\}$.

The last topic of this paper is the notion of quantum harmonic functions. 
Just as in our treatment of quantum holomorphic functions there are two
natural concepts of quantum harmonic functions called weak and strong. To
interpret the equation $\partial(\bar\partial a)=0$ with as few assumption as possible, we use Proposition \ref{aholofprop} to make the following definitions. 
An element $a\in C(D_q)$ is called strongly harmonic if $a$ is scalable and the derivative $\bar\partial a$ is antiholomorphic. Similarly $a\in C(D_q)$ is called weakly harmonic if the bilinear form $Q_{\bar\partial a}$ is equal
(on $\cD$) to a quadratic form coming from an antiholomorphic element
of $C(D_q)$.
We denote by $\text {Har}(D_q)$ the space of weakly harmonic elements of
$C(D_q)$.
It turns out that such noncommutative harmonic functions
share all the essential properties with their commutative counterparts
\cite{R}. This is summarized in the following theorem.

\begin{thm}
{\ }

\begin{enumerate}
\item An element $a\in C(D_q)$ is strongly harmonic iff $a$ is scalable and weakly harmonic.
\item  If $f\in C(D)$ is harmonic inside $D$ then the corresponding Toeplitz operator $T(f)\in C(D_q)$ is weakly harmonic.
\item  If $a\in C(D_q)$ is weakly harmonic then there exist $f\in C(D)$ which is harmonic inside $D$ such that $a=T(f)$.
\item  An element $a\in C(D_q)$ is weakly harmonic iff it can be written as
$a=a_1+a_2$ where $a_1$ is weakly holomorphic and $a_2$ is antiholomorphic.
\item (mean value) If $a=T(f)\in\text {Har}(D_q)$  then
\begin{eqnarray*}
\int_{D_q}a=\int_D f(\zeta)\,d^2\zeta=f(0).
\end{eqnarray*}
\item  (maximum principle) If $a\in\text {Har}(D_q)$  then 
$||a||_{D_q}=||\sigma(a)||_{\partial D}$.
\item  The set $\text {Har}(D_q)\subset C(D_q)$ is a closed subspace isomorphic to the Banach space $\text {Har}(D)\subset C(D)$ of continuous functions on $D$ and harmonic inside $D$.
\item  (Dirichlet problem) For every $f\in C(\partial D)$ there is a unique
weakly harmonic element $a\in C(D_q)$ such that $\sigma(a)=f$. In fact $a=T(Pf)$, where
\begin{eqnarray}
Pf(\zeta)=\int_{0}^{2\pi}\frac{1-|\zeta|^2}{ |e^{i\theta}-\zeta|^2}f(e^{i\theta})\frac{d\theta}{ 2\pi}
\end{eqnarray}
(Poisson integral)
\item  An element $a\in\text {Har}(D_q)$ is positive, $a\geq 0$, iff $a=T(f)$, where $f\in C(D)$ is harmonic inside $D$ and $f\geq 0$.
\item (Harnack's theorem) If $a_n$ is an uniformly bounded and increasing
sequence of weakly harmonic elements of $C(D_q)$, 
then $\{a_n\}$ is convergent in norm.
\end{enumerate}
\end{thm}
\proof

1. This follows from the fact that if $a$ is scalable then the quadratic
form $Q_{\bar\partial a}$ comes from the operator $\bar\partial a$.

2. From the classical harmonic analysis, if $f\in C(D)$ is harmonic inside $D$ then it is a sum $f=g+h$ where $g,h\in C(D)$ and $g$ is holomorphic and $h$ antiholomorphic inside $D$. The Toeplitz operator $T(g)$ is weakly holomorphic so $Q_{\bar\partial (T(f))}=Q_{\bar\partial (T(h))}$. But the quadratic form
$Q_{\bar\partial (T(h))}$ comes from the operator $\bar\partial (T(h))$ and, using (\ref{barpartialmono}) and (\ref{qder}), we compute:
\begin{eqnarray}
\bar\partial (T(h))=T(\delta_q(h)).\label{barpartoep}
\end{eqnarray}
But $\delta_q(h)$ is again in $C(D)$ and is antiholomorphic inside $D$.
All of this implies that the quadratic form $Q_{\bar\partial (T(f))}$
is coming from an antiholomorphic element of $C(D_q)$ as claimed.

3. From the definition, if $a\in C(D_q)$ is weakly harmonic then
the quadratic form $Q_{\bar\partial a}$ comes from an antiholomorphic
element of $C(D_q)$ which, as we know from Theorem \ref{qholothm}, is of
the form $T(k)$, where $k\in C(D)$ is antiholomorphic inside $D$.
Consider:
\begin{eqnarray*}
g(\bar\zeta):=(1-q)\bar\zeta\sum_{n=0}^\infty k(q^n\bar\zeta).
\end{eqnarray*}
It is easy to see that, just like $k$, $g\in C(D)$ and
is antiholomorphic inside $D$. Moreover it is straightforward to verify that 
$\delta_q(g)=k$. This means, using (\ref{barpartoep}) that
$\bar\partial (T(g))=T(k)$. Consequently $Q_{\bar\partial (a-T(g))}=0$,
which means that $a-T(g)$ is weakly holomorphic and as we know from Theorem \ref{qholothm}, it is of the form $T(h)$, where $h\in C(D)$ is holomorphic inside $D$. We see now that, with the above notation, $a=T(f)$, where
$f:=g+h$, and $f\in C(D)$ is harmonic inside $D$ as claimed.

4. This is a direct consequence of the proof of the previous two statements.

5. If $f\in C(D)$ is harmonic, it has the following power series expansion inside $D$:
\begin{eqnarray*}
f(\zeta)=f(0)+\sum_{n\geq 1}a_n\zeta^n + \sum_{n\geq 1}b_n\bar\zeta^n.
\end{eqnarray*}
It follows that $T(f)=f(0)I+\sum_{n\geq 1}a_nz^n + \sum_{n\geq 1}b_n\bar z^n$.
Lemma \ref{intlemma} implies that $\int_{D_q}T(f)=f(0)$, which is clearly also
the value of the integral $\int_D f(\zeta)\,d^2\zeta$.

6. From the item 4, if $a\in C(D_q)$ is weakly harmonic then it can be written as $a=a_1+a_2$ where $a_1$ is weakly holomorphic and $a_2$ is antiholomorphic. Now it is enough to apply the maximum principle to both $a_1$ and $a_2$.

7. The map $\text {Har}(D)\ni f\to T(f)\in \text {Har}(D_q)$ is linear,
one-to-one, and, by item 6, an isometry.

8. For a Toeplitz operator $T(f)$ we have $\sigma(T(f))=f|_{\partial D}$.
This, the items 2, 3, and classical results on harmonic functions \cite{R}
on the unit disk imply the thesis.

9. This just says that Toeplitz operators preserve positivity which
follows immediately from the definition.

10. This again is a direct consequence of the previous items and the classical Harnack's theorem.

$\square$

\end{document}